\documentclass[a4,12pt,draft]{article}
\usepackage{amsmath}
\usepackage{amssymb}
\usepackage{amsthm}
\newtheorem{theorem}{Theorem}[section]
\newtheorem{proposition}[theorem]{Proposition}
\newtheorem{lemma}[theorem]{Lemma}
\newtheorem{corollary}[theorem]{Corollary}
\newtheorem*{ohno}{Theorem A}
\theoremstyle{definition}
\newtheorem{remark}[theorem]{Remark}
\newtheorem*{acknowledgment}{Acknowledgment}
\begin{document}
\begin{flushleft}
\textbf{Note. (April 20, 2011.)}
\end{flushleft}
The contents of this paper are completely the same as those of 
the paper submitted to the Journal of Number Theory on June 26, 2008. 
(The current referee was chosen on August 13, 2009.) 
This paper contains many errors and ambiguities. I note that 
the papers arXiv:0908.2536v1--v7 are the author's revised versions of 
the above paper submitted to the Journal of Number Theory, i.e., 
the revised versions which do not reflect the results of the review. 
I also note that this research depends on an unpublished work of Hiroyuki Ochiai [10], 
i.e., Ochiai's proof of the sum formula for multiple zeta values. 
(Ochiai's proof can be found in, e.g., [8, pp. 60--61].) 
I talked this research at Seminar on Analytic Number Theory, 
Graduate of School of Mathematics, Nagoya University, Japan, 13 February 2008. 
The title of the talk was ``On Ohno's relation for certain multiple series, (Japanese).''
\pagebreak
\title{A parametrized generalization of Ohno's relation for multiple zeta values}
\author{Masahiro Igarashi}
\date{June 26, 2008}
\maketitle
\begin{abstract}
In this paper, we prove that certain parametrized multiple series which generalize muitiple zeta
values satisfy the same relation as Ohno's relation for multiple zeta values. 
This is a parametrized generalization of Ohno's relation for multiple zeta values. 
By virtue of this generalization, we obtain a certain equivalence between the relation for 
the parametrized multiple series and its subfamily. As applications of the above results, we obtain 
some results for multiple zeta values.
\end{abstract}
\textbf{Keywords}: Parametrized multiple series; Multiple zeta value; Ohno's relation
\section{Introduction}
The multiple zeta value (MZV, for short) is defined by the multiple series
\begin{equation*}
 \zeta(k_1,\ldots,k_n) := \sum_{0< m_1< \cdots <m_n}\frac{1}{m_1^{k_1} \cdots m_n^{k_n}},
\end{equation*}
where $k_1,\ldots ,k_n$ are positive integers and $k_n\ge2$ . Various relations for MZV's are known 
(see, e.g., \cite{S5}, \cite{S7}, \cite{S9}, \cite{S11}, \cite{S12}, \cite{S15}). 
Especially Ohno's relation for MZV's (\cite{S11}) is a well-known $\mathbb{Q}$-linear relation for MZV's.
\par
 Before we state Ohno's relation for MZV's, we recall some definitions. An index $(k_1,\ldots,k_n)$ is called 
an admissible index if it  satisfies that $k_1,\ldots,k_n$ are positive integers and $k_n \ge2$ . 
The sum $k_1+\cdots+k_n$ and $n$ are  called the weight and the depth of the index $(k_1,\ldots,k_n)$, respectively. Any admissible 
index $(k_1,\ldots,k_n)$ can be expressed as 
\begin{equation*}
\mathbf{k}: = (k_1,\ldots,k_n)
            = (\underbrace{1,\ldots,1}_{a_1},b_{1}+2,\ldots,b_{s-1}+2,\underbrace{1,\ldots,1}_{a_s},b_{s}+2),
\end{equation*}
where $a_1,\ldots,a_s,b_1,\ldots,b_s$ are non-negative integers. Then the dual index of $\mathbf{k}$ is defined by
\begin{equation*}
\mathbf{k^{'}} := ({k^{'}_1},\ldots,k^{'}_{n'})
             = (\underbrace{1,\ldots,1}_{b_s},a_{s}+2,\ldots,a_{2}+2,\underbrace{1,\ldots,1}_{b_1},a_{1}+2).
\end{equation*}
We can easily see that 
\begin{equation*}
k_1+\cdots+k_n=k^{'}_1+\cdots+k^{'}_{n'}=n+n'.
\end{equation*}
\par
Ohno's relation for MZV's is as follows.
\begin{ohno}[Ohno's relation for MZV's {\cite[Theorem 1]{S11}}]
Let\linebreak
$(k_1,\ldots,k_n)$ be an admissible index and $(k^{'}_1,\ldots,k^{'}_{n'})$ the dual index of\linebreak
$(k_1,\ldots,k_n)$. Then the identity 
\begin{equation*}
\sum_{\begin{subarray}{c}
l_1+\cdots+l_n=l\\
       l_i\in\mathbb{Z}_{\ge0}
      \end{subarray}}
 \zeta(k_1+l_1,\ldots,k_n+l_n)
=\sum_{\begin{subarray}{c}
l_1+\cdots+l_{n'}=l\\
       l_i\in\mathbb{Z}_{\ge0}
      \end{subarray}}
 \zeta(k^{'}_1+l_1,\ldots,k^{'}_{n'}+l_{n'})
\end{equation*}
holds for all integer $l\ge0$.
\end{ohno} 
Theorem A is a simultaneous generalization of the duality (\cite{S4}, \cite{S16}) 
and the sum formula for MZV's (\cite{S3}, \cite{S4}), and contains Hoffman's relation for MZV's (\cite[Theorem 5.1]{S4}). 
Therefore Theorem A yields many $\mathbb{Q}$-linear relations for MZV's. 
Some alternative proofs of Theorem A are found in \cite{S7}, \cite{S9} 
and \cite{S14} (see also Remark 2.8 in Section 2). Further D. M. Bradley (\cite[Theorem 5]{S1}) proved 
a $q$-analogue of Theorem A. This is a parametrized generalization 
of Ohno's relation for MZV's. We also prove some parametrized generalization of the relation, 
though we deal with other parametrized multiple series.
\par
Now we define the parametrized multiple series
\begin{equation*}
\begin{split}
&Z(k_1,\ldots,k_n;\alpha,\beta)\\
&:= \sum_{0\le m_1< \cdots <m_n}\frac{(\alpha)_{m_1}}{{m_1}!}\frac{{m_n}!}{(\alpha)_{m_{n}+1}}
\frac{1}{(m_1+\beta)^{k_1} \cdots (m_{n-1}+\beta)^{k_{n-1}}(m_n+\beta)^{{k_n}-1}},
\end{split}
\end{equation*}
where $(k_1,\ldots,k_n)$ is an admissible index, $\alpha,\beta\in\mathbb{C}$ with $\mathrm{Re}\,\alpha>0,\beta\notin{\mathbb{Z}_{\le0}}$,
and $(\alpha)_n$ denotes the Pochhammer symbol defined by
\begin{equation*}
(\alpha)_n = \left\{\begin{alignedat}{3}
                    \alpha(\alpha+1)&\cdots(\alpha+n-1)&\quad&\text{if $n\in\mathbb{Z}_{\ge1}$},\\
                    &1&\quad&\text{if $n=0$}.
                     \end{alignedat}
\right.
\end{equation*}
This is a generalization of MZV. Indeed, taking $\alpha=\beta=1$ in the defintion of $Z(\mathbf{k};\alpha,\beta)$, 
we get $\zeta(\mathbf{k})$. For simplicity, we denote $Z(\mathbf{k};\alpha,\alpha)$ by $Z(\mathbf{k};\alpha)$. 
In the case $\mathbf{k}=(k_1)$, we see that $Z(k_1;\alpha)=\zeta(k_1;\alpha)$, where $\zeta(s;\alpha)$ denotes 
the Hurwitz zeta-function. 
A multiple series similar to $Z(\mathbf{k};\alpha)$ was studied by M. \'{E}mery in \cite{S2}. 
His result can be regarded as a generalization of the Landen connection formula for polylogarithms.
\par
In this paper, we prove the following parametrized generalization of\linebreak
Ohno's relation for MZV's.
\begin{theorem}
\label{claim:1.1}
Let $(k_1,\ldots,k_n)$ be an admissible index and $(k^{'}_1,\ldots,k^{'}_{n'})$ the dual index of $(k_1,\ldots,k_n)$. 
Then the identity 
\begin{equation}
\label{eq:1}
\sum_{\begin{subarray}{c}l_1+\cdots+l_n=l\\
       l_i\in\mathbb{Z}_{\ge0}
      \end{subarray}}
 Z(k_1+l_1,\ldots,k_n+l_n;\alpha)
=\sum_{\begin{subarray}{c}l_1+\cdots+l_{n'}=l\\
       l_i\in\mathbb{Z}_{\ge0}
      \end{subarray}}
Z({k^{'}_1}+l_1,\ldots,{k^{'}_{n'}}+l_{n'};\alpha)
\end{equation}
holds for all integer $l\ge0$ and all complex number $\alpha$ with $\mathrm{Re}\,\alpha>0$.
\end{theorem}
Taking $\alpha = 1$ in Theorem \ref{claim:1.1}, we get Theorem A.
Therefore Theorem \ref{claim:1.1} can be regarded as a generalization of Theorem A. 
Further Theorem \ref{claim:1.1} means that $Z(\mathbf{k};\alpha)$'s satisfy many $\mathbb{Q}$-linear relations 
(see also Remark 2.10 in Section 2). 
As another example of Theorem \ref{claim:1.1}, we can obtain the following sum formula for $Z(\mathbf{k};\alpha)$'s: 
the identity 
\begin{equation*}
\zeta(n;\alpha)
=\sum_{\begin{subarray}{c}k_1+\cdots+k_{m}= n\\
       k_i\in\mathbb{Z}_{\ge1},k_m\ge2
      \end{subarray}}
Z(k_1,\ldots,k_m;\alpha)
\end{equation*}
holds for all integers $m,n$ with $0<m<n$ and all $\alpha\in\mathbb{C}$ with $\mathrm{Re}\,\alpha>0$. 
Indeed, this follows by applying (1) for the index $(k),k\in\mathbb{Z}_{\ge2}$. 
In \cite{S6}, following Ochiai's method of proving the sum formula for MZV's (\cite{S10}), 
the author proved the above sum formula (see also Remark 2.4 in Section 2). 
In this paper, we also follow Ochiai's method to prove Proposition 2.6 below, which is equivalent to Theorem \ref{claim:1.1}. 
The sum formula for MZV's was first proved by A. Granville (\cite{S3}) and D. Zagier, independently. 
Though Ochiai's proof is unpublished, it can be found in \cite{S8}. 
Generally speaking, Ochiai's method is as follows: 
first we find some multiple integral representation of generating functions of sums of multiple series; 
secondly, using the multiple integral representation, we get some duality formula for the generating functions; 
and finally we derive a relation for the multiple series from the above duality formula.
\par
By virtue of the parametrized generalization, we can obtain the following theorem.
\begin{theorem}
\label{claim:1.2}
The following assertion $(*)$ is equivalent to Theorem $\ref{claim:1.1}$: 
\par
$(*)$ Let $(k_1,\ldots,k_n)$ be an admissible index and $(k^{'}_1,\ldots,k^{'}_{n'})$ the dual index of $(k_1,\ldots,k_n)$. 
Then the identity $\mathrm{(\ref{eq:1})}$ holds for all ``even'' integer $l\ge0$ and 
all complex number $\alpha$ with $\mathrm{Re}\,\alpha>0$. 
\end{theorem}
Theorem \ref{claim:1.2} asserts an equivalence between the relation in Theorem \ref{claim:1.1} and its subfamily. 
\par
One of our motivations for studying parametrized multiple series is to apply the results to 
the study of MZV's. In this paper, as applications of the property of the parametrized multiple series $Z(\mathbf{k};\alpha)$, 
we obtain some results for MZV's (see below). 
Other applications of the property of parametrized multiple series to the study of MZV's are 
found in, e.g., \cite{S4} and \cite{S13}. 
\par
 We give an outline of the rest of this paper. 
In Section 2, we prove Theorem \ref{claim:1.1}. We first prove some properties of $Z(\mathbf{k};\alpha,\beta)$, 
especially a certain multiple integral representation (Proposition 2.2).
 Using the multiple integral representation, we prove a certain duality formula for $Z(\mathbf{k};\alpha,\beta)$'s. 
The duality formula can be regarded as a relation for generating functions 
of certain sums of $Z(\mathbf{k};\alpha)$'s (Proposition 2.5). 
Therefore we can derive a certain sum relation for $Z(\mathbf{k};\alpha)$'s (Proposition 2.6) from 
the duality formula for $Z(\mathbf{k};\alpha,\beta)$'s. 
The above method of proving Proposition 2.6 followed Ochiai's method of proving the sum formula for MZV's. 
Finally we prove the equivalence between Theorem \ref{claim:1.1} and Proposition 2.6 (Proposition 2.7). 
At the end of Section 2, as an applicaton of Theorem \ref{claim:1.1}, we get a relation for MZV's which contains 
Ohno's relation for MZV's (Corollary 2.9). 
In Section 3, using certain relations between the same sums of $Z(\mathbf{k};\alpha)$'s as in Theorem \ref{claim:1.1} 
and their derivatives (Proposition 3.2), we prove Theorem \ref{claim:1.2}. 
As an application of Theorem \ref{claim:1.2}, we obtain an equivalence between a relation for MZV's which follows from Corollary 2.9 
and its subfamily (Corollary 3.3).
\section{Proof of Theorem \ref{claim:1.1} and a relation for MZV's}
In this section, we prove Theorem \ref{claim:1.1}. First, following Ochiai's method of proving the sum formula for MZV's, 
we prove Proposition 2.6. Secondly we prove the equivalence between Theorem \ref{claim:1.1} and Proposition 2.6. 
As an application of Theorem \ref{claim:1.1}, we get a relation for MZV's which contains Ohno's relation for MZV's. 
\par
In order to prove Proposition 2.6, we prove some properties of $Z(\mathbf{k};\alpha,\beta)$.
\begin{proposition}
\label{claim:2.1}
Let $(k_1,\ldots,k_n)$ be an admissible index. Then the multiple series $Z(k_1,\ldots,k_n;\alpha,\beta)$ converges 
absolutly for $(\alpha,\beta)\in\{(\alpha,\beta)\in\mathbb{C}^2:\mathrm{Re}\,\alpha>0,\beta\notin\mathbb{Z}_{\le0}\}$ 
and uniformly in any compact subset of $\{(\alpha,\beta)\in\mathbb{C}^2:\mathrm{Re}\,\alpha>0,\beta\notin\mathbb{Z}_{\le0}\}$. 
\end{proposition}
\begin{proof}
We fix any real number $r$ with $0<r<1$ and any compact subset K of $\bigl\{\beta\in\mathbb{C}:\beta\notin\mathbb{Z}_{\le0}\bigr\}$. 
Let $(\alpha,\beta)\in\{(\alpha,\beta)\in\mathbb{C}^2:\mathrm{Re}\,\alpha{\ge}r,\beta\in\mathrm{K}\}$ and let 
$(k_1,\ldots,k_n)$ be an admissible index. For a given positive integer $m$, we first estimate the finite multiple sum 
\begin{equation}
\label{eq:2}
\sum_{0\le m_1< \cdots <m_{n-1}<m}\frac{(\alpha)_{m_1}}{{m_1}!}\frac{{m}!}{(\alpha)_{m+1}}
\frac{1}{(m_1+\beta)^{k_1} \cdots (m_{n-1}+\beta)^{k_{n-1}}(m+\beta)^{{k_n}-1}}.
\end{equation}
In the case $n=1$, we regard the sum (2) as 
\begin{equation*}
\frac{1}{(m+\alpha)(m+\beta)^{{k_1}-1}}.
\end{equation*}
Using Stirling's formula for the Gamma function, we get 
\begin{equation*}
\begin{aligned}
&\left|\sum_{0\le m_1< \cdots <m_{n-1}<m}\frac{(\alpha)_{m_1}}{{m_1}!}\frac{{m}!}{(\alpha)_{m+1}}
\frac{1}{(m_1+\beta)^{k_1} \cdots (m_{n-1}+\beta)^{k_{n-1}}(m+\beta)^{{k_n}-1}}\right|\\
&{\le}
\sum_{0\le m_1< \cdots <m_{n-1}<m}\frac{(r)_{m_1}}{{m_1}!}\frac{{m}!}{(r)_{m+1}}
\frac{1}{|m_1+\beta|^{k_1} \cdots |m_{n-1}+\beta|^{k_{n-1}}|m+\beta|^{{k_n}-1}}\\
&{\le}
\frac{{m}!}{{(r)_{m+1}}{|m+\beta|^{k_{n}-1}}}\left(\sum_{m_1=0}^{m}\frac{1}{|m_1+\beta|^{k_1}}\right)\cdots
\left(\sum_{m_{n-1}=0}^{m}\frac{1}{|m_{n-1}+\beta|^{k_{n-1}}}\right)\\
&{\ll}
\frac{1}{m^{1+r}}\left(\sum_{l=1}^{m}\frac{1}{l}\right)^{n-1}\\
&{\ll} 
\frac{(\log{m})^{n-1}}{m^{1+r}},
\end{aligned}
\end{equation*}
where the implied constants depend only on $r, (k_1,\ldots,k_n)$ and K.
Since the series 
\begin{equation*}
\sum_{m=1}^{\infty}\frac{(\log{m})^{n-1}}{m^{1+r}}
\end{equation*}
converges for $r>0$, by a theorem of Weierstrass, we get the assertion. 
\end{proof}
By Proposition \ref{claim:2.1}, we see that, for any admissible index $\mathbf{k}$, $Z(\mathbf{k};\alpha,\beta)$ is 
holomorphic in $\{(\alpha,\beta)\in\mathbb{C}^2:\mathrm{Re}\,\alpha>0,\beta\notin\mathbb{Z}_{\le0}\}$.
\par
The following multiple integral representation of $Z(\mathbf{k};\alpha,\beta)$ plays an essential role for 
the proof of Theorem \ref{claim:1.1}.
\begin{proposition}
\label{claim:2.2}
Let $(k_1,\ldots,k_n)$ be an admissible index. Then the multiple integral representation
\begin{equation}
\label{eq:3}
\begin{aligned}
&Z(k_1,\ldots,k_n;\alpha,\beta)\\
&= \idotsint\displaylimits_{1>t_{1}>\cdots>t_{k}>0}
(1-t_1)^{\alpha-1}t_{1}^{1-\beta}{\omega_{1}}(t_1)\cdots{\omega_k}(t_k)(1-t_k)^{1-\alpha}t_{k}^{\beta-1}
\,dt_1 \cdots dt_k
\end{aligned}
\end{equation}
holds for all complex numbers $\alpha,\beta$ with $\mathrm{Re}\,\alpha>0,\mathrm{Re}\,\beta>0$, 
where $k$ denotes the weight of the index $(k_1, \ldots ,k_n)$, and
\begin{equation*}
{\omega_i}(t_i) = \left\{\begin{alignedat}{3}
                    (1&-t_i)^{-1}&\quad&\text{if $\,i\in\{k_1+\cdots+k_j:j=1,\ldots,n\}$},\\
                    &{t_i}^{-1}&\quad&\text{otherwise}.
\end{alignedat}
\right.
\end{equation*}
\end{proposition}
\begin{proof}
Using the Taylor expansions at the origin of $(1-t_k)^{-\alpha}$ and $(1-t_i)^{-1}$ for $i\in\{k_1+\cdots+k_j:j=1,\ldots,n-1\}$, 
we calculate the right-hand side of (\ref{eq:3}) as follows: 
\begin{equation*}
\begin{aligned}
&\idotsint\displaylimits_{1>t_{1}>\cdots>t_{k}>0}
(1-t_1)^{\alpha-1}t_{1}^{1-\beta}{\omega_{1}}(t_1)\cdots{\omega_k}(t_k)(1-t_k)^{1-\alpha}t_{k}^{\beta-1}
\,dt_1 \cdots dt_k\\
&=\idotsint\displaylimits_{1>t_{1}>\cdots>t_{k}>0}
\frac{(1-t_1)^{\alpha-1}t_{k}^{\beta-1}\,dt_1 \cdots dt_k}{t_{1}^{\beta}t_2\cdots{t_{k_1-1}}(1-t_{k_1})t_{k_1+1}\cdots
{t_{k-1}}(1-t_k)^\alpha}\\
&=\sum_{l_1, \ldots ,\l_n\ge0}\frac{(\alpha)_{l_1}}{l_1!}
\idotsint\displaylimits_{1>t_{1}>\cdots>t_{k}>0}
\frac{(1-t_1)^{\alpha-1}t_{k_1}^{l_n}\cdots{t_{k}^{l_1+\beta-1}}\,dt_1 \cdots dt_k}
{t_{1}^{\beta}t_2\cdots{t_{k_1-1}}t_{k_1+1}\cdots
{t_{k-1}}}\\
&=\sum_{l_1, \ldots ,\l_n\ge0}
\frac{(\alpha)_{l_1}}
{l_1!}\int_{0}^{1}\frac{(1-t_1)^{\alpha-1}}{t_1^\beta}\,dt_1 
\int_{0}^{t_1}\frac{dt_2}{t_2} \cdots \int_{0}^{t_{k-2}}\frac{dt_{k-1}}{t_{k-1}}\int_{0}^{t_{k-1}}t_k^{l_1+\beta-1}\,dt_k \\
&=\sum_{l_1 ,\ldots ,\l_n\ge0}
\frac{(\alpha)_{l_1}}{l_1!}\frac{1}{(l_1+\beta)^{k_1}\cdots(l_1+\cdots+l_{n-1}+n-2+\beta)^{k_{n-1}}}\\
&\times\frac{1}{(l_1+\cdots+l_{n}+n-1+\beta)^{k_{n}-1}}\int_{0}^{1}(1-t_1)^{\alpha-1}t_1^{l_1+\cdots+l_n+n-1}\,dt_1\\
&=\sum_{l_1 ,\ldots ,\l_n\ge0}
\frac{(\alpha)_{l_1}}{l_1!}\frac{1}{(l_1+\beta)^{k_1}\cdots(l_1+\cdots+l_{n-1}+n-2+\beta)^{k_{n-1}}}\\
&\times\frac{1}{(l_1+\cdots+l_{n}+n-1+\beta)^{k_{n}-1}}\frac{\Gamma(\alpha)\Gamma(l_1+\cdots+l_n+n)}{\Gamma(\alpha+l_1+\cdots+l_n+n)}\\
&=Z(k_1,\ldots,k_n;\alpha,\beta).
\end{aligned}
\end{equation*}
The above calculation is justified by the convergence of $Z(\mathbf{k};\mathrm{Re}\,\alpha,\mathrm{Re}\,\beta)$ for an 
admissible index $(k_1,\ldots,k_n)$ and $\alpha,\beta\in\mathbb{C}$ with $\mathrm{Re}\,\alpha>0,\mathrm{Re}\,\beta>0$. 
This completes the proof of Proposition 2.2.
\end{proof}
Using (3) in Proposition 2.2, we prove the following duality for $Z(\mathbf{k};\alpha,\beta)$.
\begin{proposition}[The duality formula for $Z(\mathbf{k};\alpha,\beta)$'s]
\label{claim:2.3}
Let $\mathbf{k}$ be an admissible index and $\mathbf{k^{'}}$ the dual index of $\mathbf{k}$. Then the identity 
\begin{equation*}
Z(\mathbf{k};\alpha,\beta) = Z(\mathbf{k'};\beta,\alpha)
\end{equation*}
holds for all complex numbers $\alpha,\beta$ with $\mathrm{Re}\,\alpha>0,\mathrm{Re}\,\beta>0$. 
\end{proposition}
\begin{proof}
The proof is the same as that for MZV's in \cite[p.~510]{S16}. Indeed, the assertion follows from applying the change of 
variables 
\begin{equation*}
t_i \,\longmapsto\, 1-t_{k-i+1},
\end{equation*}
where $i=1,\ldots,k$, to the multiple integral of the right-hand side of (\ref{eq:3}) in Proposition 2.2. 
\end{proof}
Taking $\alpha = \beta$ in Proposition 2.3, we get the duality formula for $Z(\mathbf{k};\alpha)$'s.
\begin{remark}
\label{remark:1.1}
Considering the index $(k),k\in\mathbb{Z}_{\ge2}$, we can derive the following sum formula for $Z(\mathbf{k};\alpha,\beta)$'s 
from Proposition $2.3$: the identity 
\begin{equation*}
\sum_{l=0}^{\infty}\frac{1}{(l+\alpha)^m(l+\beta)^n}
=\sum_{\begin{subarray}{c}k_1+\cdots+k_{m}=m+n\\
       k_i\in\mathbb{Z}_{\ge1},k_m\ge2
      \end{subarray}}
Z(k_1,\ldots,k_m;\alpha,\beta)
\end{equation*}
holds for all integers $m, n{\ge1}$ and all $\alpha,\beta\in\mathbb{C}$ with $\mathrm{Re}\,\alpha>0,\mathrm{Re}\,\beta>0$. 
This sum formula was proved by the author ([6]) following Ochiai's method of proving the sum fomula for MZV's. 
We note that the condition $\mathrm{Re}\,\beta>0$ in the above sum formula 
can be replaced by $\beta\notin\mathbb{Z}_{\le0}$, because both sides of the above sum formula are holomorphic in 
$\mathbb{C}-\mathbb{Z}_{\le0}$ as functions of $\beta$.
\end{remark}
For simplicity, we put 
\begin{equation*}
\begin{aligned}
S_l(k_1,\ldots,k_n;\alpha) := \sum_{\begin{subarray}{c}l_1+\cdots+l_n=l\\
       l_i\in\mathbb{Z}_{\ge0}
      \end{subarray}}
 Z(k_1+l_1,\ldots,k_n+l_n;\alpha),
\end{aligned}
\end{equation*}
and define 
\begin{equation*}
\mathbf{i}_{n}^{(m)} = \left\{\begin{alignedat}{3}
                    i_{1}^{(m)}+&\cdots&+i_{n}^{(m)}\quad&\text{if $m,n\in\mathbb{Z}_{\ge1}$},\\
                    &0&\quad&\text{if $m\in\mathbb{Z}_{\ge1},n=0$},
\end{alignedat}
\right.
\end{equation*}
for $i_{1}^{(m)},\ldots,i_{n}^{(m)}\in\mathbb{Z}_{\ge0}$.
\par
By using the notation $S_l(\mathbf{k};\alpha)$, the identity (1) in Theorem \ref{claim:1.1} can be written as 
$S_l(\mathbf{k};\alpha)=S_l(\mathbf{k^{'}};\alpha)$.
\par
The following proposition asserts that $Z(\mathbf{k};\alpha,\beta)$ is a generating function of sums of $S_l(\mathbf{k};\alpha)$.
\begin{proposition}
\label{claim:2.4}
Let $(k_1,\ldots,k_n)$ be an admissible index and let $\alpha$ be a complex number with positive real part. 
Then the following two expansions hold: 
\par
$\mathrm{(i)}$
\begin{equation*}
\begin{aligned}
&Z(k_1,\ldots,k_n;\beta,\alpha)\\
&= \sum_{l=0}^{\infty}(\alpha-\beta)^l\\
&\times\sum_{i=0}^{l}\sum_{\begin{subarray}{c}i_{1}+\cdots+i_{n-1} = i\\
                                        i_{1},\ldots, i_{n-1}\in\mathbb{Z}_{\ge0}
                                       \end{subarray}}
S_{l-i}(k_1,\underbrace{1,\ldots,1}_{i_1},k_2,\ldots,k_{n-1},\underbrace{1,\ldots,1}_{i_{n-1}},k_{n};\alpha)
\end{aligned}
\end{equation*}
for all complex number $\beta$ with $|\beta-\alpha|<\mathrm{Re}\,\alpha$;
\par
$\mathrm{(ii)}$
\begin{equation*}
\begin{aligned}
&Z(k_1,\ldots,k_n;\alpha,\beta)\\
&= \sum_{l=0}^{\infty}(\alpha-\beta)^l\\
&\times\sum_{i=0}^{l}\sum_{\mathbf{i}_{k_1-1}^{(1)}+\cdots+\mathbf{i}_{k_n-2}^{(n)}=i}
S_{l-i}(k_1+\mathbf{i}_{k_1-1}^{(1)},\ldots,k_{n-1}+\mathbf{i}_{k_{n-1}-1}^{(n-1)},k_n+\mathbf{i}_{k_n-2}^{(n)};\alpha)
\end{aligned}
\end{equation*}
for all complex number $\beta$ with $|\beta-\alpha|<\mathrm{Re}\,\alpha$.
\end{proposition}
\begin{proof}
We fix any admissible index $(k_1,\ldots,k_n)$ and any complex number $\alpha_0$ with $\mathrm{Re}\,{\alpha_0}>0$. 
Then, expanding $Z(k_1,\ldots,k_n;\beta,{\alpha_0})$ into the Taylor series at $\beta={\alpha_0}$, we get 
\begin{equation}
\label{eq:4}
Z(k_1,\ldots,k_n;\beta,{\alpha_0})=\sum_{l=0}^{\infty}
\frac{1}{l!}\frac{{\mathrm{d}}^{l}}{\mathrm{d}{\beta}^{l}}
Z(k_1,\ldots,k_n;\beta,{\alpha_0})\Big|_{\beta={\alpha_0}}
(\beta-{\alpha_0})^l
\end{equation}
for all $\beta\in\mathbb{C}$ with $|\beta-{\alpha_0}|<\mathrm{Re}\,{\alpha_0}$. 
By induction on $l$, we obtain 
\begin{equation}
\label{eq:5}
\begin{aligned}
&\frac{(-1)^l}{l!}\frac{\mathrm{d}^l}{\mathrm{d}\beta^l}\frac{(\beta)_{m_1}}{(\beta)_{m_{n}+1}}\\
&=\frac{(\beta)_{m_1}}{(\beta)_{m_{n}+1}}\sum_{m_1\le{n_1}\le\cdots\le{n_l}\le{m_n}}\frac{1}{(n_1+\beta)\cdots(n_l+\beta)}\\
&=\frac{(\beta)_{m_1}}{(\beta)_{m_{n}+1}}\sum_{i=0}^{l}\sum_{i_1+\cdots+i_{n-1}=i}\sum_{\begin{subarray}{c}l_1+\cdots+l_n\\
                                                          +j_{1}^{(1)}+\cdots+j_{i_1}^{(1)}\\
                                                        \cdots\\
+j_{1}^{(n-1)}+\cdots+j_{i_{n-1}}^{(n-1)}=l-i
\end{subarray}}
\frac{1}{(m_1+\beta)^{l_1}\cdots(m_n+\beta)^{l_n}}\\
&\times\sum_{\begin{subarray}{c}
m_1<m_{11}<\cdots<m_{1 i_1}<m_2\\
\cdots\\
m_{n-1}<m_{n-1 1}<\cdots<m_{n-1 i_{n-1}}<m_n
\end{subarray}}
\prod_{p=1}^{n-1}\prod_{q=1}^{i_p}\frac{1}{(m_{pq}+\beta)^{j_{q}^{(p)}+1}}
\end{aligned}
\end{equation}
for all integer $l\ge1$.
Using (\ref{eq:4}) and (\ref{eq:5}), we get (i). 
\par
Similarly, expanding $Z(k_1,\ldots,k_n;{\alpha_0},\beta)$ into the Taylor series at $\beta={\alpha_0}$, we get 
\begin{equation}
\label{eq:6}
Z(k_1,\ldots,k_n;{\alpha_0},\beta)=\sum_{l=0}^{\infty}
\frac{1}{l!}\frac{{\mathrm{d}}^{l}}{\mathrm{d}{\beta}^{l}}
Z(k_1,\ldots,k_n;{\alpha_0},\beta)\Big|_{\beta={\alpha_0}}
(\beta-{\alpha_0})^l
\end{equation}
for all $\beta\in\mathbb{C}$ with $|\beta-{\alpha_0}|<\mathrm{Re}\,{\alpha_0}$. 
We can easily verify that 
\begin{equation}
\label{eq:7}
\begin{aligned}
&\frac{(-1)^l}{l!}\frac{\mathrm{d}^l}{\mathrm{d}\beta^l}\frac{1}{(m_1+\beta)^{k_1}\cdots(m_{n-1}+\beta)^{k_{n-1}}(m_n+\beta)^{k_{n}-1}}\\
&=\sum_{i=0}^{l}
\sum_{\mathbf{i}_{k_1-1}^{(1)}+\cdots+\mathbf{i}_{k_n-2}^{(n)}=i }
\sum_{l_1+\cdots+l_n=l-i}\frac{1}{(m_n+\beta)^{k_n-1+\mathbf{i}_{k_n-2}^{(n)}+l_{n}}}\\
&\times\prod_{j=1}^{n-1}\frac{1}{(m_j+\beta)^{k_j+\mathbf{i}_{k_j-1}^{(j)}+l_j}}
\end{aligned}
\end{equation}
holds for any admissible index $(k_1,\ldots,k_n)$ and all integer $l\ge0$.
Using (\ref{eq:6}) and (\ref{eq:7}), we get (ii). 
\end{proof}
Using the above properties of $Z(\mathbf{k};\alpha,\beta)$, we prove the following proposition.
\begin{proposition}
\label{claim:2.5}
Let $(k_1,\ldots,k_n)$ be an admissible index and $(k^{'}_1,\ldots,k^{'}_{n'})$ the dual index of $(k_1,\ldots,k_n)$. 
Then the identity 
\begin{equation}
\begin{aligned}
\label{eq:8}
&\sum_{i=0}^{l}\sum_{\begin{subarray}{c}i_{1}+\cdots+i_{n-1} =i\\
                                        i_{1},\ldots, i_{n-1}\in\mathbb{Z}_{\ge0}
                                       \end{subarray}}
S_{l-i}(k_1,\underbrace{1,\ldots,1}_{i_1},k_2,\ldots,k_{n-1},\underbrace{1,\ldots,1}_{i_{n-1}},k_{n};\alpha)\\ 
&=\sum_{i=0}^{l}\sum_{\mathbf{i}_{k^{'}_1-1}^{(1)}+\cdots+\mathbf{i}_{k^{'}_{n'}-2}^{(n')}=i}
S_{l-i}({k^{'}_1}+\mathbf{i}_{{k^{'}_1}-1}^{(1)},\ldots,{k^{'}_{n'-1}}+\mathbf{i}_{k^{'}_{n'-1}-1}^{(n'-1)},
{k^{'}_{n'}}+\mathbf{i}_{k^{'}_{n'}-2}^{(n')};\alpha)
\end{aligned}
\end{equation}
holds for all integer $l\ge0$ and all complex number $\alpha$ with $\mathrm{Re}\,\alpha>0$.
\end{proposition}
\begin{proof}
By applying Proposition \ref{claim:2.3} to Proposition \ref{claim:2.4}, the generating functions of both sides of 
(\ref{eq:8}) coincide. Therefore we get the assertion. 
\end{proof}
Theorem \ref{claim:1.1} follows from Proposition \ref{claim:2.5}. In fact, these are equivalent. 
\begin{proposition}
\label{claim:2.6}
Theorem $\ref{claim:1.1}$ and Propositon $\ref{claim:2.5}$ are equivalent.
\end{proposition}
\begin{proof}
Let $(k_1,\ldots,k_n)$ be an admissible index, $(k^{'}_1,\ldots,k^{'}_{n'})$ the dual index of $(k_1,\ldots,k_n)$, 
and $i_1,\ldots,i_{n-1}$ non-negative integers. 
Then we note that, by the definition of dual index, the dual index of 
\begin{equation*}
(k_1,\underbrace{1,\ldots,1}_{i_1},k_2,\ldots,k_{n-1},\underbrace{1,\ldots,1}_{i_{n-1}},k_n)
\end{equation*}
takes the form of 
\begin{equation*}
({k^{'}_1}+\mathbf{i}_{{k^{'}_1}-1}^{(1)},\ldots,{k^{'}_{n'-1}}+\mathbf{i}_{{k^{'}_{n'-1}}-1}^{(n'-1)},
{k^{'}_{n'}}+\mathbf{i}_{{k^{'}_{n'}-2}}^{(n')}).
\end{equation*}
Further the identity 
\begin{equation*}
i_1+\cdots+i_{n-1}=\mathbf{i}_{{k'_1}-1}^{(1)}+\cdots+\mathbf{i}_{{k^{'}_{n'-1}}-1}^{(n'-1)}+\mathbf{i}_{{k'_{n'}-2}}^{(n')}
\end{equation*}
holds, because an admissible index and its dual index have the same weight. 
By these facts, we can easily prove that Theorem \ref{claim:1.1} implies Proposition \ref{claim:2.5}. 
\par
Conversly we suppose that Proposition \ref{claim:2.5} is true. Then, by the above facts, we can rewrite (\ref{eq:8}) as 
\begin{equation}
\label{eq:9}
\begin{aligned}
&S_l(k_1,\ldots,k_n;\alpha)-S_l({k^{'}_1},\ldots,{k^{'}_{n'}};\alpha)\\
&=-\sum_{i=1}^{l}\sum_{\begin{subarray}{c}i_{1}+\cdots+i_{n-1} = i\\
                                        i_{1},\ldots, i_{n-1}\in\mathbb{Z}_{\ge0}
                                       \end{subarray}}
\Bigl\{S_{l-i}(\mathbf{k}_{i_1,\ldots,i_{n-1}};\alpha)-S_{l-i}(\mathbf{k}^{'}_{i_1,\ldots,i_{n-1}};\alpha)\Bigr\},
\end{aligned}
\end{equation}
where 
\begin{equation*}
\mathbf{k}_{i_1,\ldots,i_{n-1}}:=(k_1,\underbrace{1,\ldots,1}_{i_1},k_2,\ldots,k_{n-1},\underbrace{1,\ldots,1}_{i_{n-1}},k_n),
\end{equation*}
and $\mathbf{k}^{'}_{i_1,\ldots,i_{n-1}}$ is the dual index of $\mathbf{k}_{i_1,\ldots,i_{n-1}}$.
Proposition \ref{claim:2.5} contains the duality formula for $Z(\mathbf{k};\alpha)$'s. 
Therefore, by using (\ref{eq:9}) and induction on $l$, we get Theorem \ref{claim:1.1}. 
This completes the proof of Proposition 2.7. 
\end{proof}
\begin{remark}
Taking $\alpha=1$ in Proposition $\mathrm{2.6}$, we get a relation for MZV's. 
By using the relation and the same argument as in the proof of Proposition $\mathrm{2.7}$, 
we can obtain Ohno's relation for MZV's. This is an alternative proof of Ohno's relation for MZV's.
\end{remark}
By Theorem \ref{claim:1.1} and expanding ${S_l}(\mathbf{k};\alpha)$ into the Taylor series, we get the following.
\begin{corollary}
Let $\alpha_0$ be a complex number with positive real part, $\mathbf{k}$ an admissible index, and $\mathbf{k^{'}}$ 
the dual index of $\mathbf{k}$. 
Then the identity 
\begin{equation*}
\frac{\mathrm{d}^{m}}{\mathrm{d}\alpha^{m}}S_{l}(\mathbf{k};\alpha){\Bigl|}_{\alpha={{\alpha}_0}}
=\frac{\mathrm{d}^{m}}{\mathrm{d}\alpha^{m}}S_{l}(\mathbf{k^{'}};\alpha){\Bigl|}_{\alpha={{\alpha}_0}}
\end{equation*}
holds for all integers $l,m\ge0$.
\end{corollary}
Using (\ref{eq:5}) and (\ref{eq:7}) in the proof of Proposition 2.5, we obtain 
\begin{equation*}
\begin{aligned}
&\frac{(-1)^m}{m!}\frac{\mathrm{d}^{m}}{\mathrm{d}\alpha^{m}}S_{l}(k_1,\ldots,k_n;\alpha)\\
&=\sum_{\begin{subarray}{c}l_{1}+\cdots+l_{n} = l\\
                            l_{i}\in\mathbb{Z}_{\ge0}
                                   \end{subarray}}
\sum_{\begin{subarray}{c}i_{1}+\cdots+i_{n} \\
                         +{\mathbf{i}}_{k_1+l_1}^{(1)}+\cdots+{\mathbf{i}}_{k_n+l_{n}-1}^{(n)}=m
                                   \end{subarray}}
S_{i_n}(k_{1}+l_{n}+{\mathbf{i}}_{k_1+l_1}^{(1)},\underbrace{1,\ldots,1}_{i_1},k_{2}+l_{2}+{\mathbf{i}}_{k_2+l_2}^{(2)},\\
&\ldots,k_{n-1}+l_{n-1}+{\mathbf{i}}_{k_{n-1}+l_{n-1}}^{(n-1)},\underbrace{1,\ldots,1}_{i_{n-1}},
k_{n}+l_{n}+{\mathbf{i}}_{k_n+l_{n}-1}^{(n)};\alpha)
\end{aligned}
\end{equation*}
for all integers $l,m{\ge0}$. Therefore the identity in Corollary 2.9 is a relation for $Z(\mathbf{k};{\alpha}_0)$'s. 
\par
Taking ${\alpha}_{0}=1$ in Corollary 2.9, we get a relation for MZV's. In \cite{S9}, G. Kawashima proved a relation 
for MZV's which contains Ohno's relation for MZV's. Our relation for MZV's also contains Ohno's relation for MZV's.
\begin{remark}
It is well-known that $\zeta(4)=4\zeta(1,3)$. However the functions $Z(4;\alpha)(=\zeta(4;\alpha))$ and $4Z(1,3;\alpha)$ 
are not identically equal in the half-plane $\{\alpha\in\mathbb{C}:\mathrm{Re}\,\alpha>0\}$. 
Indeed, we can easily verify that $Z(4;2)\neq4Z(1,3;2)$. Similarly $\zeta(4)=4\zeta(2,2)/3$, but 
$Z(4;\alpha)\neq4Z(2,2;\alpha)/3$ for some $\alpha\in\mathbb{C}$ with $\mathrm{Re}\,\alpha>0$. 
These differences probably come from the absence of suitable expressions of the product of two $Z(\mathbf{k};\alpha)$'s for $\alpha\neq1$. 
In the case $\alpha=1$, see \cite{S7}.
\end{remark}
\section{Proof of Theorem \ref{claim:1.2} and an equivalence of a relation for MZV's}
In this section, using Proposition 3.2 below, we prove Theorem \ref{claim:1.2}. As an application of Theorem \ref{claim:1.2}, we prove 
an equivalence for the relation in Corollary 2.9. 
\par
For each complex number $\beta$ with $\mathrm{Re}\,\beta>0$, we put 
\begin{equation*}
\mathrm{D}(\beta):=\{\alpha\in\mathbb{C}:|\alpha-\beta|<\mathrm{Re}\,\beta/2\}.
\end{equation*}
\par
We first prove a lemma.
\begin{lemma}
Let $(k_1,\ldots,k_{n})$ be an admissible index and let $\beta$ be a complex number with positive real part. 
Then the series
\begin{equation}
\label{eq:10}
\sum_{l=0}^{\infty}(\alpha-\beta)^l
\sum_{i=0}^{l}\sum_{\begin{subarray}{c}i_{1}+\cdots+i_{n-1} = i\\
                                        i_{1},\ldots, i_{n-1}\in\mathbb{Z}_{\ge0}
                                       \end{subarray}}
S_{l-i}(k_1,\underbrace{1,\ldots,1}_{i_1},k_2,\ldots,k_{n-1},\underbrace{1,\ldots,1}_{i_{n-1}},k_{n};\alpha) 
\end{equation}
and
\begin{equation}
\label{eq:11}
\sum_{l=0}^{\infty}(\alpha-\beta)^l
\sum_{i=0}^{l}\sum_{\mathbf{i}_{k_1-1}^{(1)}+\cdots+\mathbf{i}_{k_n-2}^{(n)}=i}
S_{l-i}(k_1+\mathbf{i}_{k_1-1}^{(1)},\ldots,k_{n-1}+\mathbf{i}_{k_{n-1}-1}^{(n-1)},k_n+\mathbf{i}_{k_n-2}^{(n)};\alpha)
\end{equation}
converge absolutly for $\alpha\in\mathrm{D(\beta)}$ and uniformly in any compact subset of $\mathrm{D(\beta)}$.
\end{lemma}
\begin{proof}
For any fixed real number $r$ with $\mathrm{Re}\,\beta/2<r<{3\mathrm{Re}\,\beta}/2$, we put 
\begin{equation*}
\mathrm{D}_r(\beta):=\mathrm{D}(\beta)\cap\{\alpha\in\mathbb{C}:\mathrm{Re}\,\alpha{\ge}r\}.
\end{equation*}
Then we get 
\begin{equation*}
\begin{aligned}
&\left|\sum_{i=0}^{l}\sum_{\begin{subarray}{c}i_{1}+\cdots+i_{n-1} = i\\
                                        i_{1},\ldots, i_{n-1}\in\mathbb{Z}_{\ge0}
                                       \end{subarray}}
S_{l-i}(k_1,\underbrace{1,\ldots,1}_{i_1},k_2,\ldots,k_{n-1},\underbrace{1,\ldots,1}_{i_{n-1}},k_{n};\alpha)\right|\\
&\le\sum_{i=0}^{l}\sum_{\begin{subarray}{c}i_{1}+\cdots+i_{n-1} = i\\
                                        i_{1},\ldots, i_{n-1}\in\mathbb{Z}_{\ge0}
                                       \end{subarray}}
S_{l-i}(k_1,\underbrace{1,\ldots,1}_{i_1},k_2,\ldots,k_{n-1},\underbrace{1,\ldots,1}_{i_{n-1}},k_{n};r)\\
&=\frac{(-1)^l}{l!}\frac{\mathrm{d}^l}{\mathrm{d}z^l}Z(k_1,\ldots,k_n;z,r){\Bigl|}_{z=r}
\end{aligned}
\end{equation*}
for all $\alpha\in{\mathrm{D}_r(\beta)}$ and all integer $l\ge0$. 
Further, by Cauchy's theorem, we get 
\begin{equation*}
\begin{aligned}
&\frac{(-1)^l}{l!}\frac{\mathrm{d}^l}{\mathrm{d}z^l}Z(k_1,\ldots,k_n;z,r){\Bigl|}_{z=r}\\
&=\frac{(-1)^l}{2{\pi}{\sqrt{-1}}}\int_{|z-r|=\rho}\frac{Z(k_1,\ldots,k_n;z,r)}{(z-r)^{l+1}}\,dz\\
&\le\frac{Z(k_1,\ldots,k_n;r-\rho,r)}{\rho^l},
\end{aligned}
\end{equation*}
where $\rho\in\mathbb{R}$ with $\mathrm{Re}\,\beta/2<\rho<r$.
Thus we get 
\begin{equation*}
\begin{aligned}
&\left|\alpha-\beta\right|^l\left|\sum_{i=0}^{l}\sum_{\begin{subarray}{c}i_{1}+\cdots+i_{n-1} = i\\
                                        i_{1},\ldots, i_{n-1}\in\mathbb{Z}_{\ge0}
                                       \end{subarray}}
S_{l-i}(k_1,\underbrace{1,\ldots,1}_{i_1},k_2,\ldots,k_{n-1},\underbrace{1,\ldots,1}_{i_{n-1}},k_{n};\alpha)\right|\\
&<\left(\frac{\mathrm{Re}\,\beta}{2\rho}\right)^lZ(k_1,\ldots,k_n;r-\rho,r)
\end{aligned}
\end{equation*}
for all $\alpha{\in}\mathrm{D}_r(\beta)$, all integer $l\ge0$ and 
a fixed $\rho\in\mathbb{R}$ with $\mathrm{Re}\,\beta/2<\rho<r$. 
By the above estimate and a theorem of Weierstrass, we get the assertion for (\ref{eq:10}). 
\par
By the same argument as above, we can prove the assertion for (\ref{eq:11}).
\end{proof}
By Lemma 3.1, we see that (\ref{eq:10}) and (\ref{eq:11}) are holomorphic in $\mathrm{D}(\beta)$ as functions of $\alpha$.
\par
We shall use the following proposition to prove Theorem 1.2.
\begin{proposition}
Let $(k_1,\ldots,k_n)$ be an admissible index and let $\beta$ be a complex number with positive real part. 
Then, for all integer $m\ge0$, the following two identities hold: 
\begin{align}
\nonumber
&\sum_{i=0}^{m}\sum_{\mathbf{i}_{k_1-1}^{(1)}+\cdots+\mathbf{i}_{k_n-2}^{(n)}=i}
S_{m-i}(k_1+\mathbf{i}_{k_1-1}^{(1)},\ldots,k_{n-1}+\mathbf{i}_{k_{n-1}-1}^{(n-1)},k_n+\mathbf{i}_{k_n-2}^{(n)};\beta)\\
\label{eq:12}
&=(-1)^{m}\sum_{l=0}^{m}\frac{1}{(m-l)!}\frac{\mathrm{d}^{m-l}}{\mathrm{d}\alpha^{m-l}}\Biggl\{
\sum_{i=0}^{l}
\sum_{\begin{subarray}{c}i_{1}+\cdots+i_{n-1} =i\\
                                        i_{1},\ldots, i_{n-1}\in\mathbb{Z}_{\ge0}
                                       \end{subarray}}
S_{l-i}(k_1,\underbrace{1,\ldots,1}_{i_1},k_2,\\
\nonumber
&\ldots,k_{n-1},\underbrace{1,\ldots,1}_{i_{n-1}},k_{n};\alpha)\Biggr\}
{\Biggr|}_{\alpha=\beta}
\end{align}
and
\begin{align}
\nonumber
&\sum_{i=0}^{m}\sum_{\begin{subarray}{c}i_{1}+\cdots+i_{n-1} = i\\
                                        i_{1},\ldots, i_{n-1}\in\mathbb{Z}_{\ge0}
                                       \end{subarray}}
S_{m-i}(k_1,\underbrace{1,\ldots,1}_{i_1},k_2,\ldots,k_{n-1},\underbrace{1,\ldots,1}_{i_{n-1}},k_{n};\beta)\\
\label{eq:13}
&=(-1)^{m}\sum_{l=0}^{m}\frac{1}{(m-l)!}\frac{\mathrm{d}^{m-l}}{\mathrm{d}\alpha^{m-l}}\Biggl\{
\sum_{i=0}^{l}
\sum_{\mathbf{i}_{k_1-1}^{(1)}+\cdots+\mathbf{i}_{k_n-2}^{(n)}=i}
S_{l-i}(k_1+\mathbf{i}_{k_1-1}^{(1)},\\
\nonumber
&\ldots,k_{n-1}+\mathbf{i}_{k_{n-1}-1}^{(n-1)},k_n+\mathbf{i}_{k_n-2}^{(n)};\alpha)\Biggr\}
{\Biggr|}_{\alpha=\beta}.
\end{align}
\end{proposition}
\begin{proof}
We fix any admissible index $(k_1,\ldots,k_n)$ and any complex number $\beta_0$ with $\mathrm{Re}\,{\beta_0}>0$. 
Then, by Proposition 2.5 (i), the expansion
\begin{equation}
\label{eq:14}
\begin{aligned}
&Z(k_1,\ldots,k_n;{\beta_0},\alpha)\\
&=\sum_{l=0}^{\infty}(\alpha-{\beta_0})^l\\
&\times\sum_{i=0}^{l}\sum_{\begin{subarray}{c}i_{1}+\cdots+i_{n-1} = i\\
                                        i_{1},\ldots, i_{n-1}\in\mathbb{Z}_{\ge0}
                                       \end{subarray}}
S_{l-i}(k_1,\underbrace{1,\ldots,1}_{i_1},k_2,\ldots,k_{n-1},\underbrace{1,\ldots,1}_{i_{n-1}},k_{n};\alpha)
\end{aligned}
\end{equation}
holds for all $\alpha\in\mathrm{D}(\beta_0)$. Further, by Lemma 3.1, we can differentiate the right-hand 
side of (\ref{eq:14}) term by term with respect to $\alpha$ in $\mathrm{D}({\beta}_{0})$. Thus we obtain 
\begin{equation*}
\begin{aligned}
&\frac{(-1)^m}{m!}\frac{\mathrm{\mathrm{d}}^m}{\mathrm{\mathrm{d}}\alpha^m}Z(k_1,\ldots,k_n;{\beta}_{0},\alpha){\Bigl|}_{\alpha={{\beta}_{0}}}\\
&=(-1)^{m}\sum_{l=0}^{m}\frac{1}{(m-l)!}\frac{\mathrm{d}^{m-l}}{\mathrm{d}\alpha^{m-l}}\Biggl\{
\sum_{i=0}^{l}
\sum_{\begin{subarray}{c}i_{1}+\cdots+i_{n-1} = i\\
                                        i_{1},\ldots, i_{n-1}\in\mathbb{Z}_{\ge0}
                                       \end{subarray}}
S_{l-i}(k_1,\underbrace{1,\ldots,1}_{i_1},k_2,\\
&\ldots,k_{n-1},\underbrace{1,\ldots,1}_{i_{n-1}},k_{n};\alpha)\Biggr\}{\Biggl|}_{\alpha={{\beta}_{0}}}
\end{aligned}
\end{equation*}
for all integer $m\ge0$. In the proof of Proposition 2.5, we proved that the left-hand side of the above identity 
is equal to that of (\ref{eq:12}). This completes the proof of (\ref{eq:12}). 
\par
By Proposition 2.5 (ii), Lemma 3.1, and the same argument as above, we can prove (\ref{eq:13}).
\end{proof}
Now we prove Theorem \ref{claim:1.2}. 
\begin{proof}[Proof of Theorem $\ref{claim:1.2}$]
It is trivial that Theorem \ref{claim:1.1} implies Theorem \ref{claim:1.2} $(*)$.
\par
Conversly we suppose that Theorem \ref{claim:1.2} $(*)$ 
is true. Then it is enough to prove that, for all positive ``odd'' integer $m$ and all $\alpha\in\mathbb{C}$ with 
$\mathrm{Re}\,\alpha>0$, the identity 
\begin{equation}
\label{eq:15}
S_{m}(\mathbf{k};\alpha)=S_{m}(\mathbf{k}^{'};\alpha)
\end{equation}
holds, where $\mathbf{k}$ is any admissible index, and $\mathbf{k}^{'}$ is the dual index of $\mathbf{k}$. 
\par
Using the identities (\ref{eq:12}) for $\mathbf{k}$ = $(k_1,\ldots,k_n)$ and (\ref{eq:13}) for 
$\mathbf{k}^{'}$= $(k^{'}_1,\ldots,k^{'}_{n'})$, and recalling what we noted the form of 
the dual index of a certain index in the proof of Propositon 2.7, we get 
\begin{align}
\nonumber
&\sum_{i=0}^{m}\sum_{\begin{subarray}{c}i_{1}+\cdots+i_{n'-1} = i\\
                                        i_{1},\ldots, i_{n'-1}\in\mathbb{Z}_{\ge0}
                                       \end{subarray}}
\Bigl\{S_{m-i}(\mathbf{h}^{'}_{i_1,\ldots,i_{n'-1}};\beta)-S_{m-i}(\mathbf{h}_{i_1,\ldots,i_{n'-1}};\beta)\Bigr\}\\
\label{eq:16}
&=(-1)^{m}\sum_{l=0}^{m}\frac{1}{(m-l)!}\\
\nonumber
&\times\frac{\mathrm{d}^{m-l}}{\mathrm{d}\alpha^{m-l}}\Biggl[\sum_{i=0}^{l}\sum_{\begin{subarray}{c}i_{1}+\cdots+i_{n-1} =i\\
                                        i_{1},\ldots, i_{n-1}\in\mathbb{Z}_{\ge0}
                                       \end{subarray}}
\Bigl\{S_{l-i}(\mathbf{k}_{i_1,\ldots,i_{n-1}};\alpha)-S_{l-i}(\mathbf{k}^{'}_{i_1,\ldots,i_{n-1}};\alpha)\Bigr\}{\Biggr]}{\Biggl|}_{\alpha=\beta} 
\end{align}
for all $\beta\in\mathbb{C}$ with $\mathrm{Re}\,\beta>0$ and all integer $m\ge0$, 
where $\mathbf{k}_{i_1,\ldots,i_{n-1}}$ and $\mathbf{k}^{'}_{i_1,\ldots,i_{n-1}}$ are the same as in the proof of Proposition 2.7,
\begin{equation*}
\mathbf{h}_{i_1,\ldots,i_{n'-1}}:=(k^{'}_1,\underbrace{1,\ldots,1}_{i_1},k^{'}_2,\ldots,k^{'}_{n'-1},\underbrace{1,\ldots,1}_{i_{n'-1}},
k^{'}_{n'}),
\end{equation*}
and $\mathbf{h}^{'}_{i_1,\ldots,i_{n'-1}}$ is the dual index of $\mathbf{h}_{i_1,\ldots,i_{n'-1}}$. 
The identity (\ref{eq:16}) can be rewritten as 
\begin{align}
\nonumber
&\Bigl\{1+(-1)^{m+1}\Bigr\}\Bigl\{S_m(\mathbf{k};\beta)-S_m(\mathbf{k'};\beta)\Bigr\}\\
\nonumber
&=-\sum_{i=1}^{m}\sum_{\begin{subarray}{c}i_{1}+\cdots+i_{n'-1} =i\\
                                        i_{1},\ldots, i_{n'-1}\in\mathbb{Z}_{\ge0}
                                       \end{subarray}}
\Bigl\{S_{m-i}(\mathbf{h}^{'}_{i_1,\ldots,i_{n'-1}};\beta)-S_{m-i}(\mathbf{h}_{i_1,\ldots,i_{n'-1}};\beta)\Bigr\}\\
\label{eq:17}
&+(-1)^{m}\sum_{i=1}^{m}\sum_{\begin{subarray}{c}i_{1}+\cdots+i_{n-1} =i\\
                                        i_{1},\ldots, i_{n-1}\in\mathbb{Z}_{\ge0}
                                       \end{subarray}}
\Bigl\{S_{m-i}(\mathbf{k}_{i_1,\ldots,i_{n-1}};\beta)-S_{m-i}(\mathbf{k}^{'}_{i_1,\ldots,i_{n-1}};\beta)\Bigr\}\\
\nonumber
&+(-1)^{m}\sum_{l=0}^{m-1}\frac{1}{(m-l)!}\\
\nonumber
&\times\frac{\mathrm{d}^{m-l}}{\mathrm{d}\alpha^{m-l}}\Biggl[\sum_{i=0}^{l}\sum_{\begin{subarray}{c}i_{1}+\cdots+i_{n-1} =i\\
                                        i_{1},\ldots, i_{n-1}\in\mathbb{Z}_{\ge0}
                                       \end{subarray}}
\Bigl\{S_{l-i}(\mathbf{k}_{i_1,\ldots,i_{n-1}};\alpha)-S_{l-i}(\mathbf{k}^{'}_{i_1,\ldots,i_{n-1}};\alpha)\Bigr\}{\Biggr]}{\Biggl|}_{\alpha=\beta}. 
\end{align}
 Further we suppose that $m$ is a positive odd integer. Then, applying 
Theorem \ref{claim:1.2} $(*)$ to the right-hand 
side of (\ref{eq:17}), 
we get 
\begin{align}
\nonumber
&2\Bigl\{S_m(\mathbf{k};\beta)-S_m(\mathbf{k'};\beta)\Bigr\}\\
\nonumber
&=-\sum_{\begin{subarray}{c}i=0\\
i:\text{odd}
\end{subarray}}^{m-1}\sum_{\begin{subarray}{c}i_{1}+\cdots+i_{n'-1} =m-i\\
                                        i_{1},\ldots, i_{n'-1}\in\mathbb{Z}_{\ge0}
                                       \end{subarray}}
\Bigl\{S_{i}(\mathbf{h}^{'}_{i_1,\ldots,i_{n'-1}};\beta)-S_{i}(\mathbf{h}_{i_1,\ldots,i_{n'-1}};\beta)\Bigr\}\\
\label{eq:18}
&-\sum_{\begin{subarray}{c}i=0\\
i:\text{odd}
\end{subarray}}^{m-1}\sum_{\begin{subarray}{c}i_{1}+\cdots+i_{n-1} =m-i\\
                                        i_{1},\ldots, i_{n-1}\in\mathbb{Z}_{\ge0}
                                       \end{subarray}}
\Bigl\{S_{i}(\mathbf{k}_{i_1,\ldots,i_{n-1}};\beta)-S_{i}(\mathbf{k}^{'}_{i_1,\ldots,i_{n-1}};\beta)\Bigr\}\\
\nonumber
&-\sum_{l=0}^{m-1}\frac{1}{(m-l)!}\\
\nonumber
&\times\frac{\mathrm{d}^{m-l}}{\mathrm{d}\alpha^{m-l}}\Biggl[\sum_{\begin{subarray}{c}i=0\\
i:\text{odd}
\end{subarray}}^{l}\sum_{\begin{subarray}{c}i_{1}+\cdots+i_{n-1} =l-i\\
                                        i_{1},\ldots, i_{n-1}\in\mathbb{Z}_{\ge0}
                                       \end{subarray}}
\Bigl\{S_{i}(\mathbf{k}_{i_1,\ldots,i_{n-1}};\alpha)-S_{i}(\mathbf{k}^{'}_{i_1,\ldots,i_{n-1}};\alpha)\Bigr\}{\Biggr]}{\Biggl|}_{\alpha=\beta}. 
\end{align}
Therefore, by using (\ref{eq:18}) and induction on the positive odd integer $m$, we get (\ref{eq:15}). 
This completes the proof of Theorem \ref{claim:1.2}. 
\end{proof}
As an application of Theorem \ref{claim:1.2}, we prove the following equivalence for the relation in Corollary 2.9.
\begin{corollary}
Let $\alpha_0$ and $\alpha_1$ be complex numbers with positive real parts. 
Then the following two assertions are equivalent:
\par
$\mathrm{(i)}$
Let $\mathbf{k}$ be an admissible index and $\mathbf{k^{'}}$ the dual index of $\mathbf{k}$. Then the identity
\begin{equation*}
\frac{\mathrm{d}^{m}}{\mathrm{d}\alpha^{m}}S_{l}(\mathbf{k};\alpha){\Bigl|}_{\alpha={{\alpha}_0}}
=\frac{\mathrm{d}^{m}}{\mathrm{d}\alpha^{m}}S_{l}(\mathbf{k^{'}};\alpha){\Bigl|}_{\alpha={{\alpha}_0}}
\end{equation*}
holds for all integers $l,m\ge0$. 
\par
$\mathrm{(ii)}$
Let $\mathbf{k}$ be an admissible index and $\mathbf{k^{'}}$ the dual index of $\mathbf{k}$. Then the identity
\begin{equation*}
\frac{\mathrm{d}^{m}}{\mathrm{d}\alpha^{m}}S_{l}(\mathbf{k};\alpha){\Bigl|}_{\alpha={{\alpha}_1}}
=\frac{\mathrm{d}^{m}}{\mathrm{d}\alpha^{m}}S_{l}(\mathbf{k^{'}};\alpha){\Bigl|}_{\alpha={{\alpha}_1}}
\end{equation*}
holds for all ``even'' integer $l\ge0$ and all integer $m\ge0$.
\end{corollary}
\begin{proof}
We suppose that (i) is true. 
Then, expanding $S_{l}(\mathbf{k};\alpha)$ into the Taylor series at $\alpha={\alpha}_{0}$, 
we see that, for all integer $l\ge0$ and all $\alpha\in\mathbb{C}$ with $|\alpha-{{\alpha}_{0}}|<\mathrm{Re}\,{{\alpha}_{0}}$, 
the identity 
\begin{equation*}
S_{l}(\mathbf{k};\alpha)=S_{l}(\mathbf{k^{'}};\alpha)
\end{equation*}
holds, where $\mathbf{k}$ is an admissible index, and $\mathbf{k^{'}}$ is the dual index of $\mathbf{k}$. 
By the uniqueness theorem for analytic functions, the above identity holds for all $\alpha\in\mathbb{C}$ with $\mathrm{Re}\,{\alpha}>0$. 
This is exactly Theorem \ref{claim:1.1}. Therefore, by Theorem \ref{claim:1.2}, we get Theorem \ref{claim:1.2} $(*)$. 
Clearly Theorem \ref{claim:1.2} $(*)$ implies (ii). 
\par
By the same argument as above, we can prove that (ii) implies (i). 
\end{proof}
Taking ${{\alpha}_0}={{\alpha}_1}=1$ in Corollary 3.3, we get an equivalence between the relation for MZV's 
which we stated at the end of Section 2 and its subfamily.
\begin{acknowledgment}
The author would like to express his gratitude to Professor Yoshio Tanigawa for his useful advice.
\end{acknowledgment}
\begin{flushleft}
Graduate School of Mathematics\\
Nagoya University\\
Furo-cho, Chikusa-ku, Nagoya 464-8602, Japan\\
\textbf{E-mail address}: m05003x@math.nagoya-u.ac.jp\\
\end{flushleft}
\end{document}